\documentclass{fic-l}
\usepackage{amsmath}
\usepackage{amsthm}
\usepackage{amsfonts}
\usepackage{enumerate}

\newcommand{\fa}{\mbox{$\mathfrak{a}$}}

\newcommand{\ord}{\mathop{\rm ord}}

\newcommand{\Z}{\mathop{\bf Z}}

\newcommand{\cyc}[1]{\mathop{{\bf Q}(\zeta_{#1})}}

\newcommand{\tw}[1]{\mathop{\theta_{#1}}}
\newcommand{\twx}[1]{\mathop{\theta_{#1}(x)}}

\newcommand{\twb}[1]{\mathop{\theta_{[{#1}]}}}
\newcommand{\tz}[1]{\mathop{\theta(\zeta_{{#1}}})}
\newcommand{\tzbar}[1]{\mathop{\theta(\overline{\zeta_{{#1}}}})}
\newcommand{\twbx}[1]{\mathop{\theta_{[{#1}]}(x)}}

\newtheorem{theorem}{Theorem}[section]

\title{On the existence of cyclic difference sets with small parameters}
\author{Leonard D. Baumert}
\address{325 Acero Place\\
Arroyo Grande, CA 93420}

\author{Daniel M. Gordon}
\address{IDA Center for Communications Research\\
4320 Westerra Court\\
San Diego, CA 92121}
\email{gordon@ccrwest.org}

\subjclass{Primary 05B10}

\dedicatory{This paper is dedicated to Hugh Williams on the occasion
  of his 60th birthday.}

\begin{document}

\begin{abstract}
  Previous surveys by Baumert \cite{baumert} and Lopez and Sanchez
  \cite{ls} have resolved the existence of cyclic $(v,k,\lambda)$
  difference sets with $k \leq 150$, except for six open cases.  In this
  paper we show that four of those difference sets do not exist.  We
  also look at the existence of difference sets with $k\leq 300$, and
  cyclic Hadamard difference sets with $v \leq 10{,}000$.  Finally, we
  extend \cite{gordon} to show that no cyclic projective planes exist
  with non-prime power orders $\leq 2 \cdot 10^9$.
\end{abstract}

\maketitle

\section{Introduction}

A $(v,k,\lambda)$ difference set is a subset
$D =\{d_1,d_2,\ldots,d_k\}$ of a group $G$
such that each nonidentity element $g\in G$ can be represented 
as $g=d_i d_j^{-1}$ in
exactly $\lambda$ ways.  In this paper we will be concerned with
{\em cyclic} difference sets, where $G$ will be taken to be the cyclic group
$\Z / v \Z$.  The {\em order} of a difference set is $n=k-\lambda$.

Baumert \cite{baumert} gave a complete list of parameters for
cyclic difference sets
with $k \leq 100$.  
Lander gave a table of possible abelian difference set parameters
with $k \leq50$.  
Kopilovich \cite{kopilovich} extended the search to $k<100$, and 
Lopez and Sanchez \cite{ls} looked at all possible parameters for
abelian difference sets with 
$k\leq 150$.  
Table~\ref{tab:ds150} shows their
open cases for cyclic difference sets, four of which we show do
not exist.


\begin{table}[htbp]
  \centering
\begin{tabular}[c]{|cccccc|}
\hline
$v$ & $k$ & $\lambda$ & $n$ & Status & Reference\\ \hline
$429$ & $108$ & $27$& $81$ &  No & Theorem~\ref{thm:429}\\
$715$ & $120$ & $20$& $100$ &  No & Theorem 4.20 of \cite{lander} \\
$351$ & $126$ & $45$& $81$ & No & Schmidt Test \cite{schmidt2}\\
$837$ & $133$ & $21$  & $112$  & No &  Schmidt Test \cite{schmidt2} \\
$419$ & $133$ & $42$  & $91$  & Open &  \\
$465$ & $145$ & $45$ & $100$  & Open &  \\
\hline
\end{tabular}
  \caption{Possible Difference Sets with $k\leq 150$}
  \label{tab:ds150}
\end{table}

In addition to settling some of these open cases, 
we have extended these calculations to larger values of $k$, using the
same procedure of applying the numerous known necessary conditions.
The open cases for $k \leq 300$ are given in Tables~\ref{tab:ds300a}
and~\ref{tab:ds300b}.  The cases with $\gcd(v,n)$ greater
than one are given separately, because of Ryser's conjecture that no
cyclic difference sets exist with $\gcd(v,n)>1$.


\begin{table}[htbp]
  \centering
\begin{tabular}[c]{|ccccc|}
\hline
$v$ & $k$ & $\lambda$ & $n$ & $\gcd(v,n)$ \\ \hline
$945$ & $177$ & $33$ & $144$ & $9$   \\
$5859$ & $203$ & $7$ & $196$ & $7$   \\
$1785$ & $224$ & $28$ & $196$ & $7$   \\
$2574$ & $249$ & $24$ & $225$ & $9$   \\
$2160$ & $255$ & $30$ & $225$ & $45$   \\
$1925$ & $260$ & $35$ & $225$ & $25$   \\
\hline
\end{tabular}
  \caption{Possible CDS with 
$150 \leq k\leq 300$ and $\gcd(v,n)>1$}
  \label{tab:ds300a}
\end{table}


\begin{table}[htbp]
  \centering
\begin{tabular}[c]{|cccc|}
\hline
$v$ & $k$ & $\lambda$ & $n$  \\ \hline
$1123$ & $154$ & $21$ & $133$     \\
$645$ & $161$ & $40$ & $121$  \\
$1093$ & $169$ & $26$ & $143$     \\
$1111$ & $186$ & $31$ & $155$     \\
$469$ & $208$ & $92$ & $116$     \\
$1801$ & $225$ & $28$ & $197$     \\
$2291$ & $230$ & $23$ & $207$     \\
$639$ & $232$ & $84$ & $148$     \\
$2869$ & $240$ & $20$ & $220$     \\
$1381$ & $276$ & $55$ & $221$     \\
$817$ & $289$ & $102$ & $187$     \\
$781$ & $300$ & $115$ & $185$     \\
\hline
\end{tabular}
  \caption{Possible CDS with 
$150 \leq k\leq 300$ and $\gcd(v,n)=1$}
  \label{tab:ds300b}
\end{table}

We will give the details of computations that excluded possible
difference sets in these tables.  Most of the techniques are well
known, and are described briefly in Section~\ref{sec:nec}.  
A few parameters require more effort, such as the $(429,108,27)$
difference set which is shown not to exist in Section~\ref{sec:429}.

In Section~\ref{sec:hadamard} we look at cyclic Hadamard
difference sets, with
$v=4n-1$, $k=2n-1$, $\lambda=n-1$.
There are three known families, and it is
conjectured that no others exist. 

In Section~\ref{sec:ppc} we look at the Prime Power Conjecture, which
states that all abelian difference sets with $\lambda=1$ have $n$ a
prime power.  For the cyclic case we extend earlier computations by
the second author \cite{gordon} to $2 \cdot 10^9$, showing that no
such difference sets exist when $n$ is not a prime power.

Details of the computations, such as nonexistence proofs for the hard
cases of cyclic projective planes mentioned in Section~\ref{sec:ppc},
are not included in the paper.
A web site
{\tt http://www.ccrwest.org/diffsets.html}, maintained by the second
author, lists many
known difference sets and gives nonexistence proofs.

\section{Necessary Conditions}  \label{sec:nec}

As in other searches (\cite{baumert69}, \cite{kopilovich},
\cite{lander}, \cite{ls}) we will go through values of 
$(v,k,\lambda)$ up to a given $k$, applying known necessary
conditions to eliminate most parameters, and dealing with survivors on
a case-by-case basis. 
By a simple counting argument
we must have $(v-1)\lambda = k(k-1)$.  We may assume $k \leq v/2$,
since the complement of a $(v,k,\lambda)$ difference set is a 
$(v,v-k,v-2k+\lambda)$ difference set.  Some other conditions (see
\cite{jungnickel} for references) are:

\begin{theorem}\label{thm:testA}
(Schutzenberger) If $v$ is even, $n$ must be a square.
\end{theorem}

\begin{theorem}
(Bruck-Chowla-Ryser) If $v$ is odd, the equation
$$
nX^2 + (-1)^{(v-1)/2} \lambda Y^2 = Z^2
$$
must have a nontrivial integer solution.
\end{theorem}

\begin{theorem}
(Mann) 
If $w>1$ is a divisor of $v$, $p$ is a prime divisor of $n$, 
$p^2$ does not divide $n$, and $p^j \equiv -1 \pmod w$, then
no $(v,k,\lambda)$ difference set exists.
\end{theorem}

\begin{theorem}
(Arasu \cite{arasu})
If $w>1$ is a divisor of $v$, $p$ is a prime divisor of $n$, 
and
\begin{itemize}
\item $\gcd(v,k)=1$,
\item $n$ is a nonsquare,
\item $\gcd(p,v)=1$,
\item $p$ is a multiplier,
\end{itemize}
then 
$$
w v (-1)^{(v/w-1)/2}
$$
is a square in the ring of $p$-adic integers.
\end{theorem}

Baumert gives four necessary tests used for his search in
\cite{baumert69}, which include Theorem~\ref{thm:testA} and 
three theorems of Yamamoto \cite{yamamoto}.
Lander gives a number of conditions in Chapter 4 of \cite{lander}.
Ones that are used to exclude possible difference sets
include Theorems 4.19, 4.20, 4.27,
4.30, 4.31, 4.32, 4.33, and 4.38.

\section{Constructing the Tables}\label{sec:tables}

To extend previous tables of possible cyclic difference sets, we apply
the theorems of the previous section to eliminate most possible
parameters.  Ones that survive these tests are dealt with on a case by
case basis.  In this section we give methods from \cite{baumert} for
dealing with certain difficult cases, and show an example of their
application.

\subsection{Polynomial Congruences}  \label{sec:baumert}

Let $\theta(x)$ be the difference set polynomial
$$
\theta(x) = x^{d_1} + x^{d_2} + \ldots + x^{d_k},
$$
and $\zeta_{v}$ be a primitive $v$th root of unity.
Then $D$ is a difference set if and only if
$$
\tz{v} \tzbar{v} = n.
$$

In \cite{baumert69} and 
\cite{baumert} a method is given for constructing or showing the
nonexistence of difference sets.  Define
$$
\tw{w} (x) \equiv \theta(x) \pmod {f_w(x)}
$$
where $f_w(x)$ is the $w$th cyclotomic polynomial, and
$$
\twb{w} (x) \equiv \theta(x) \pmod {x^w-1}.
$$
The method is based on the congruences proved in \cite{baumert}:
\begin{equation}\label{eq:rums1}
w \twb{w} \equiv w \tw{w} - 
\sum_{\begin{array}{c}r|w\\r \ne w\end{array}} \mu(w/r)
r \left(\twb{r} - \tw{w}\right) \frac{x^w-1}{x^r-1} \pmod {x^w-1},
\end{equation}
and
\begin{equation}\label{eq:rums2}
\tw{w} \equiv \twb{w/p} \pmod {p,f_{w_1}^{p^{a-1}}}
\end{equation}
where $w = p^a w_1$, with $\gcd(p,w_1)=1$.

Thus, given $\tw{w}$ for a divisor of $v$ and $\twb{r}$ for all 
divisors $r$ of $w$, one may compute $\twb{w}$.  To find $\tw{w}$, we 
may use the equation
$$
\tz{w} \tzbar{w} = n.
$$
Furthermore, if $\fa$ is an ideal in $\cyc{w}$ 
for which 
$\fa \overline{\fa} = (n)$
with generator $\sum a_i \zeta_w^i$, 
then if 
$\theta_{w}(\zeta_w) \in \fa$ we have
\begin{equation}
  \label{eq:kron}
\twx{w} = \pm x^j \sum a_i x^i  
\end{equation}
by a theorem of Kronecker that any algebraic integer, all of whose
conjugates have absolute value 1, must be a root of unity.

So to determine the existence of a particular $(v,k,\lambda)$
difference set, we may factor $n$ in cyclotomic fields $\cyc{w}$ for
$w|v$, and apply congruences (\ref{eq:rums1}) and (\ref{eq:rums2}) to
construct $\twb{w}$ or show that none exists.  This approach was used
by Howard Rumsey to prove the nonexistence of difference sets
$(441,56,7)$ and $(891,90,9)$ \cite{baumert69}.  


\subsection{Contracted Multipliers}
The following two theorems, both proved in 
\cite{baumert}, are very useful:

Let 
$$\twbx{w} = b_0 + b_1x + \ldots b_{w-1}x^{w-1}.
$$
The following is Lemma 3.8 of \cite{baumert}:

\begin{theorem}\label{thm:baumert}
  For every divisor $w$ of $v$, there exists integers 
$b_i \in [0,v/w]$ such that
\begin{equation}
  \label{eq:lin}
\sum_{i=0}^{w-1} b_i = k,
\end{equation}
\begin{equation}
  \label{eq:sq}
\sum_{i=0}^{w-1} b_i^2 = n+\lambda v/w,
\end{equation}
and
\begin{equation}
  \label{eq:sq2}
\sum_{i=0}^{w-1} b_i b_{i-j}  = \lambda v/w
\end{equation}
for $j=1,\ldots,w-1$, where $i-j$ is taken modulo $w$.
\end{theorem}

The $b_i$'s are the number of $d_j$'s in $D$ satisfying 
$d_j \equiv i \pmod w$.  
These equations often are sufficient to show nonexistence of a 
difference set.  When they are not, we may sometimes use
multipliers to get further conditions.

A $w$-multiplier of a difference set is an integer $t$ prime to $w$ for 
which there is an integer $s$ such that
$$
\theta(x^t) \equiv x^s \theta(x) \pmod {x^w-1}.
$$
The following is Theorem 3.2 in \cite{baumert}, and a
generalization is given as Theorem 5.6 in \cite{lander}.

\begin{theorem}\label{thm:wmult}
Let $D$ be a $(v,k,\lambda)$ cyclic difference set
with
$n=p_1^{\alpha_1} p_2^{\alpha_2} \cdots
p_s^{\alpha_s}$.  Let $w$ be a divisor of $v$
and $t$ be an integer relatively prime to $w$.  
If for $i=1,2,\ldots,s$ there is an
integer $j=j(i)$ such that
$$
p_i^{j} \equiv t \pmod w,
$$
then $t$ is a $w$-multiplier of $D$.
\end{theorem}

If we have a $w$-multiplier for $D$, this gives us further
restrictions on the $b_i$'s, since if $i$ and $j$ are in the
same orbit of $t$ modulo $w$, then we must have $b_i = b_j$.

\subsection{Using Contracted Multipliers}\label{sec:429}

As an example of using these methods to eliminate a possible cyclic
difference set, consider the 
first open case of Ryser's conjecture, $(429,108,27)$.  

\begin{theorem}\label{thm:429}
  No $(429,108,27)$ difference set exists.
\end{theorem}

\begin{proof}
By Theorem~\ref{thm:wmult}, $3$ is a 143-multiplier.

The orbits of the residues modulo 143 have sizes
$$
1^1     3^4     5^2     15^8.
$$

Let
$$
\theta_{[143]}(x) = c_0 + c_1 x + \ldots c_{142} x^{142}.
$$

From Theorem~\ref{thm:baumert} we have $\sum c_i = k = 108$, and 
$\sum c_i^2 = n+\lambda v/w = 162$, so
$$
  162 = c_0^2 + 15(c_1^2 +\ldots + c_{29}^2) + 
5(c_{13}^2 +  c_{26}^2) 
+ 3(c_{11}^2 +c_{22}^2 + c_{44}^2 + c_{77}^2) 
$$
and
$$
  \sum c_i c_{i+j} = 81, \ \ \ {\rm for}\  j=1,\ldots,142
$$

There are $14,896$ solutions to the first equation, and 
a quick computer
search shows that none of these satisfy the second.
  
\end{proof}


This method still works when $w=v$.
For example, consider a $(303,151,75)$ difference set.  By
Theorem~\ref{thm:wmult}, which for $w=v$ is known as 
the Second Multiplier Theorem, 16 is a multiplier, with three orbits
of size 1 and 12 orbits of size 25.  Therefore a difference set would
have to be a union of one of the size-1 orbits and six of the size-25
ones.  None of these 2772 possibilities form a difference set, and so
no (303,151,75) difference set exists.
Several other similar cases are given in 
Table~\ref{tab:elim}.


\begin{table}[htbp]
  \centering
\begin{tabular}[c]{|cccccc|}
\hline
$v$ & $k$ & $\lambda$ &multiplier  & $w$ & solutions to (\ref{eq:lin})
and (\ref{eq:sq})\\ \hline
$429$ & $ 108$ & $27$ & 3 & 143 & 14896 \\
303 & 151 & 75 & 16 & 303 & 2772 \\
2585 & 153 & 9 & 2 & 235 & 0 \\
$616$ & $165$ & $44$ &11 & 56 & 301485532\\
407 & 175 & 75 & 2 & 37 & 0 \\
4401 & 176 & 7 & 13 & 489 & 504 \\
544 & 181 & 60 & 3 & 68 & 96 \\
3949 & 189 & 9 & 3 & 3949 & 2 \\
1545 & 193 & 24 & 8 & 515 & 0 \\
$1380$ & $197$ & $28$ & 2 & 115 & 0  \\
1609 & 201 & 25 & 2 & 1609 & 8 \\
6271 & 210 & 7 & 29 & 6271 & 30 \\
1056 & 211 & 42 & 13 & 44 & 6240 \\
2233 & 217 & 21 & 16 & 319 & 8512 \\
6301 & 225 & 8 & 31 & 6301 & 0 \\
601 & 225 & 84 & 3 & 601 & 56 \\
595 & 243 & 99 & 2 & 119 & 216 \\
611 & 245 & 98 & 2 & 47 & 0 \\
2057 & 257 & 32 & 3 & 187 & 0 \\
2591 & 260 & 26 & 3 & 2591 & 10 \\
3181 & 265 & 22 & 3 & 3181 & 12 \\
1061 & 265 & 66 & 199 & 1061 & 4 \\
531 & 265 & 132 & 4 & 177 & 0 \\
1615 & 270 & 45 & 4 & 323 & 17024 \\
2691 & 270 & 27 & 3 & 299 & 114592 \\
28325 & 292 & 3 & 2 & 103 & 0 \\
591 & 295 & 147 & 16 & 591 & 2772 \\
10990 & 297 & 8 & 9 & 157 & 0 \\
\hline
\end{tabular}
  \caption{Cases eliminated by Theorem~\ref{thm:baumert}}
  \label{tab:elim}
\end{table}

\subsection{Schmidt's Test}  \label{sec:schmidt}

Schmidt (\cite{schmidt}, \cite{schmidt2}) has shown that, under certain
conditions, 
a root of unity times
$\theta(\zeta_v)$ must be in a subfield of $\cyc{v}$. 
For a prime $q$ and integer $m$ with prime factorization
$\prod_{i=1}^t p_i^{c_i}$, let
$$
m_q = 
\left \{
\begin{array}{ll}
\prod_{p_i \ne q} p_i & \mbox{if $m$ is odd or $q=2$,}\\
4 \prod_{p_i \ne 2,q} p_i & \mbox{otherwise.}
\end{array}
\right.
$$

Define $F(m,n)=\prod_{i=1}^t p_i^{b_i}$ to be the minimum multiple of
$\prod_{i=1}^t p_i$ such that for every pair $(i,q)$, $i \in
\{1,\ldots,t\}$, $q$ a prime divisor of $n$, at least one of the
following conditions is satisfied:
\begin{enumerate}
\item $q=p_i$ and $(p_i,b_i) \neq (2,1)$,
\item $b_i = c_i$,
\item $q \neq p_i$ and $q^{\ord_{m_q}(q)} \not \equiv 1 \pmod {p_i^{b_i+1}}$.
\end{enumerate}

Schmidt then shows

\begin{theorem}
Assume $|X|^2=n$ for $X \in \Z[\zeta_v]$.  Then
$X \zeta_v^j \in \Z[\zeta_{F(v,n)}]$ for some $j$.
\end{theorem}

When $F(v,n)$ is significantly less than $v$, this theorem gives a powerful
condition on the difference set.  Schmidt uses it to show

\begin{theorem}\label{thm:schmidt}
  For a $(v,k,\lambda)$ cyclic difference set, we have
$$
n \leq \frac{F(v,n)^2}{4 \varphi(F(v,n))},
$$
where $\varphi$ denote's Euler's totient function.  
\end{theorem}

Theorem~\ref{thm:schmidt} eliminates 29 difference sets with $k \leq 300$.

Very recently, Leung, Ma and Schmidt \cite{lms} have shown that no
cyclic difference set exists with order $n$ a power of a prime $>3$
and
$(n,v)>1$.
This eliminates the difference set $(505,225,100)$.  They also
eliminate certain cases for powers of 3, such as $(2691,270,27)$.

\section{Cyclic Hadamard Difference Sets}\label{sec:hadamard}

A cyclic Hadamard difference set is a difference set with parameters
$v=4n-1$, $k=2n-1$, $\lambda=n-1$.  All known cyclic Hadamard
difference sets are of one of the following types:
\begin{enumerate}
\item $v$ prime.
\item $v$ a product of twin primes.
\item $v = 2^n-1$.
\end{enumerate}

It has been conjectured that no others exist.  
Song and Golomb \cite{song_golomb} excluded all but 17 cases up 
to $v=10,000$.  Kim and Song \cite{kim_song} eliminated four of those.  
The remaining ones
are listed in Table~\ref{tab:chds}, along with 
their current status.  Six can be shown not to exist by theorems in
Lander's book~\cite{lander}.


\begin{table}[htbp]
  \centering
\begin{tabular}[c]{|ccccc|}
\hline
$v$ & $k$ & $\lambda$ & Status & Comment\\ \hline
$3439$ & $1719$ & $859$ & Open &  \\
$4355$ & $2177$ & $1088$ & Open &  \\
$4623$ & $2311$ & $1155$ & No & Thm. 4.19 of \cite{lander} \\
$5775$ & $2887$ & $1443$ & No &  Thm. 4.19 of \cite{lander} \\
$7395$ & $3697$ & $1848$ & No &  Thm. 4.20 of \cite{lander} \\
$7743$ & $3871$ & $1935$ & No &  Thm. 4.19 of \cite{lander} \\
$8227$ & $4113$ & $2056$ & No &  Thm. 4.20 of \cite{lander}\\
$8463$ & $4231$ & $2115$ & No &  Thm. 4.19 of \cite{lander} \\
$8591$ & $4295$ & $2147$ & Open &  \\
$8835$ & $4417$ & $2208$ & Open &  \\
$9135$ & $4567$ & $2283$ & Open &  \\
$9215$ & $4607$ & $2303$ & Open &  \\
$9423$ & $4711$ & $2355$ & Open &  \\
\hline
\end{tabular}
  \caption{Open Cases  for Cyclic Hadamard Difference Sets}
  \label{tab:chds}
\end{table}

\section{Cyclic Projective Planes} \label{sec:ppc}

A difference set with $\lambda=1$ is called a planar difference set.  
The Prime Power Conjecture (PPC) states that all abelian planar difference
sets have order $n$ a prime power.  In \cite{gordon}, it was shown
that the PPC is true for $n < 2{,}000{,}000$.

Since that paper, several developments have made it possible to extend
those computations.  Faster computers with more memory are part of it,
but also 64-bit computing allow calcuations to be done in
single-precision, which results in a large speedup.  Using the methods
of \cite{gordon}, we have shown that no cyclic planar difference sets
of non-prime power order $n$ exist with $n < 2 \cdot 10^9$.

Most orders can be eliminated by various quick tests given in
\cite{gordon}.  There were $605$ orders which survived these tests,
and were dealt with using a theorem of Evans and Mann \cite{em}
(Lander \cite{lander} proved a generalization for abelian groups):

\begin{theorem}\label{thm:mann}
Let $D$ be a $(v,k,1)$ planar cyclic difference set of order 
$n = k-1$.  If
$t_1$, $t_2$, $t_3$, and $t_4$ are numerical multipliers such that
$$
t_1-t_2 \equiv t_3-t_4 \pmod {v},
$$
then $v$ divides the least common multiple of $(t_1-t_2,t_1-t_3)$.
\end{theorem}

In \cite{gordon} this theorem was used to create a hash table for
differences $t_i - t_j$ less than one million, to find a collision
that could be used to eliminate an order.
For orders up to $2 \cdot 10^9$, all but two could be eliminated with
differences up to $4 \cdot 10^8$.  The two most difficult were 
$n=40027523$ and $n=883007071$.  These were finally eliminated with
pairs with differences $420511455$ and $164204313$, respectively.

\bibliography{diff}
\bibliographystyle{plain}

\end{document}